\documentclass[]{ipart}

\startlocaldefs
\theoremstyle{plain}
\newtheorem{thm}{Theorem}[section]
\newtheorem{definition}{Definition}[section]
\newtheorem{remark}{Remark}[section]

\def\R{\mathbb{R}}
\def\S{\Sigma}
\def\p{\partial}
\def\tr{\textmd{tr}}
\newcommand{\m}{\mathfrak{m}}
\newcommand{\ric}{\mathrm{Ric}}
\def\og{\mathring \gamma}
\def\onabla{\mathring \nabla}
\endlocaldefs

\begin{document}

\begin{frontmatter}

\title{An Overview of Bartnik's Quasi-Local Mass}

\begin{aug}
    \author{\fnms{Stephen} \snm{McCormick}\ead[label=e1]{stephen.m.mccormick@ltu.se}}
\address{Department of Engineering Sciences and Mathematics\\
		Luleå University of Technology\\
		971 87 Luleå\\
	Sweden\\
	\printead{e1}}
\end{aug}

\begin{center}
	\textit{Dedicated to the memory of Robert Bartnik} 
\end{center}

\vspace{-5mm}
\begin{abstract}
	This article provides a concise introduction to Bartnik's quasi-local mass, and surveys a selection of results pertaining to the understanding of it. The aim is to serve as both an entry point to the topic, and a quick reference of results for those already familiar with it.

\end{abstract}



\tableofcontents

\end{frontmatter}

\section{Introduction}
In 1982, Roger Penrose published a list of unsolved problems in classical general relativity \cite{penrose1982}, emphasising that important mathematical works in classical general relativity -- as opposed to some more speculative theories -- would have a permanent place in physics. The first of fourteen on this list of open problems is to \textit{``find a suitable quasi-local definition of energy-momentum in general relativity''}. The problem remains open and under active investigation, with a few candidates still under consideration for a good definition. The subject of this article one such definition, namely that given by Robert Bartnik \cite{Bartnik-89}. The purpose of this article is not only to concisely survey a selection of results on the Bartnik mass, but also to provide an introduction or entry point to the topic.

To understand the problem of quasi-local mass, we begin with two basic facts of general relativity. First, the total mass of an asymptotically flat initial data set satisfying the dominant energy condition is well-defined, and given by a surface integral at infinity known as the ADM mass. Second, there cannot exist a local mass density for the gravitational field, which is essentially a consequence of the equivalence principle. The problem of quasi-local mass can be viewed as a compromise in a sense, asking for a suitable definition of the mass within a bounded domain $\Omega$ in an initial data set, rather than at a point.

In order to be a suitable measure of mass in general relativity, it should of course correspond to the physical notion of mass in some sense. However, understanding what this means is highly non-trivial. However, motivated by this, one may expect that a good quasi-local mass definition should satisfy certain properties. For example, it is generally agreed upon that a good definition should at the very least satisfy:

\begin{enumerate}
	\item[(I)] \textit{Positivity}. The mass should always be non-negative.
	\item[(II)] \textit{Compatibility with the ADM mass}. The limit of the mass over an appropriate exhaustion of an asymptotically flat manifold should yield the ADM mass.
\end{enumerate}

There are also other natural criteria that one may insist that a good definition of quasi-local mass should satisfy. However, the community is not in full agreement of exactly what these should be. Some of the most common conditions, and the ones we will look more closely at, are as follows.
\begin{enumerate}
	\item[(III)] \textit{Monotonicity}. Given $\Omega_1\subset\Omega_2$, the mass of $\Omega_1$ should be no more than the mass of $\Omega_2$.
	\item[(IV)] \textit{Small sphere limit}. The limit of the mass taken over small spheres, divided by the volume of a small sphere, should recover the mass density of the matter. There are also expectations on the next order terms, which we will elaborate on in Section \ref{sec:criteria}.
	\item[(V)] \textit{Rigidity}. If the mass of some domain satisfying the dominant energy condition is zero, then it should be flat space.
\end{enumerate}

We will reserve discussion of how the Bartnik mass measures up to these properties until Section \ref{sec:criteria}, after some background and definitions are given.

Over the years there has been an abundance of different definitions of quasi-local mass, with most eventually falling out of favour for a variety of reasons. While this article focuses only on the definition of Bartnik, the reader is directed to the ``living review'' of Szabados \cite{LRR} for a comprehensive survey of the different quasi-local mass definitions.

We give a precise definition of the Bartnik mass in Section \ref{sec:defn}, however Bartnik's idea is essentially the following. Let $\Omega\subset M$ be a bounded domain in some initial data set $(M,\gamma,K)$, and consider a set $\widehat{\mathcal A}$ of asymptotically flat initial data sets $(M_o,\gamma_o,K_o)$ each containing a copy of $(\Omega,\gamma_\Omega,K_\Omega)$, defined in an appropriate sense. The Bartnik mass is taken to be the infimum of the ADM mass over all such $(M_o,\gamma_o,K_o)$ satisfying physically motivated conditions. The idea is that the Bartnik mass localises the ADM mass, in the sense that it is the ADM mass of an initial data set whose only contribution to the total mass comes from $\Omega$.

As the Bartnik mass is given in terms of this infimum, there is no direct way to compute it. This is a major challenge for the Bartnik mass that has driven a good deal of research since it was first given, as we only can compute the value in a few special cases. When Bartnik first posed this definition, he conjectured that the infimum would be realised by initial data that is stationary in the region exterior to $\Omega$. There is a great deal of evidence to support this conjecture in a broad sense, and it is even further conjectured that this stationary extension is unique. This gives hope that the Bartnik mass can be computed by solving for this unique extension and finding the mass of that. However, in the broadest sense this conjecture is now known to be false due to work of Anderson and Jauregui \cite{A-J}, which we will discuss in more detail below.

Note that an initial data set realising the infimum is not expected to be smooth across $\partial \Omega$ in general, and the best we can expect is Lipschitz here. This motivates a modification to the definition of Bartnik mass to allow the extensions of a domain $\Omega$ to share this low regularity there. That is, to take the infimum over a larger set of manifolds that are not required to be smooth across $\partial \Omega$. By doing so, one can then consider the stationary extension problem as a boundary value problem, which we discuss in more detail in Section \ref{sec:statstat} (see, in particular, Section \ref{sec:existunique}). Throughout this article we will avoid any discussion of the regularity required for the results presented and all objects should be considered smooth, unless otherwise specified.

\section{Asymptotically flat manifolds and the ADM mass}
For the sake of completeness, we give some background on the ADM mass. Note that the dimension $n$ is taken to be no less than $3$ throughout.

\begin{definition} \label{def:AF}
	A Riemannian $n$-manifold $(M,\gamma)$ (possibly with boundary) is said to be \textit{asymptotically flat} if there is a compact set $K\subset M$ and a diffeomorphism $\varphi:\mathbb{R}^n\setminus \overline B\to M\setminus K$, where $\overline B$ is a closed ball, satisfying
	\begin{align*}
		\varphi^*\gamma-\og=O(|x|^{-\tau}), \qquad \onabla (\varphi^*\gamma)=O(|x|^{-\tau-1}),\qquad \onabla^2 (\varphi^*\gamma)=O(|x|^{-\tau-2}),
	\end{align*}
	where $\og$ and $\onabla$ are the standard metric and connection on $\mathbb R^n$, and $\tau>\frac{n-2}2$.
\end{definition}
Note that our definition of an asymptotically flat manifold only permits one asymptotic end, that is, one copy of $\varphi:\mathbb{R}^n\setminus \overline B$. This is purely for simplicity and convenience.

On an asymptotically flat manifold, we can define a quantity known as the ADM mass, which is viewed as the total mass of a system from the perspective of general relativity.
\begin{definition}
	The ADM mass $\m_{ADM}$ of an asymptotically flat manifold $(M,\gamma)$ is given by
	\begin{equation}\label{eq:admmassdefn}
		\m_{ADM}(M,\gamma)=\frac{1}{2(n-1)\omega_{n-1}}\lim_{R\to\infty}\int_{S_R}\,\onabla^i (\varphi^*\gamma_{ik})-\onabla_k (\og^{ij}\varphi^*\gamma_{ij})\,\nu^k\,dS,
	\end{equation} 
	where $S_R\subset \mathbb R^n\setminus \overline B$ is a sphere of radius $R$, $\nu$ is the outer unit normal and $\omega_{n-1}$ is the volume of the unit $(n-1)$-sphere.
\end{definition}
By the work of Bartnik \cite{Bartnik-86} and Chru\'sciel \cite{chrusciel1986mass} it is now well-known that the ADM mass is well-defined and independent of the diffeomorphism $\varphi$ used to define it. For this reason we simply write, by an abuse of notation, $\gamma$ to denote the original metric and its pullback under such a diffeomorphism. Note that asymptotic flatness alone is not sufficient to ensure that the ADM mass is in fact finite, however with the additional assumption that scalar curvature is integrable it turns out to be the case.

In fact, the ADM mass is closely related to scalar curvature, so we briefly note some properties. First note that the linearisation of the scalar curvature operator at $\gamma$ acting on a symmetric two-tensor $h$, is given by
\begin{equation}\label{eq:DR}
	D_\gamma R(h)=\nabla^i\nabla^j(h_{ij})-\Delta(\gamma^{ij}h_{ij})-\ric^{ij}h_{ij}.
\end{equation}
We can use this expression, the divergence theorem and a Taylor expansion about $\og$ to write the ADM mass as
\begin{equation}\label{eq:scalarmass}
	\m_{ADM}(M,\gamma)=\frac{1}{2(n-1)\omega_{n-1}}\int_{M\setminus B_{R_0}}R(\gamma)\,dV_{\og}+C,
\end{equation}
where $B_{R_0}$ is some fixed large ball and $C$ denotes some collection of terms that are finite on an asymptotically flat manifold. This is an easy way to see that the ADM mass is finite if and only if the scalar curvature is in $L^1$. One can also obtain \eqref{eq:scalarmass} by expressing scalar curvature in terms of the metric and its derivatives, and noting that the highest order term is a divergence whose integral agrees with the ADM mass. We will come back to the linearisation of scalar curvature \eqref{eq:DR} in Section \ref{sec:statstat}.

We now state the well-known and celebrated positive mass theorem (see Remark \ref{rem:pmt} below for some historical comments).
\begin{thm}\label{thm:PMT1}
	Let $(M,\gamma)$ be a complete asymptotically flat manifold (without boundary) satisfying $R(\gamma)\geq0$ and $R(\gamma)\in L^1(M)$ then $\m_{ADM}(M,\gamma)\geq0$ with equality if and only if $(M,\gamma)$ is isometric to Euclidean space.
\end{thm}
\begin{remark}
	An asymptotically flat manifold may be viewed as initial data for general relativity that is \textit{time-symmetric}. In this case, the scalar curvature is proportional to the local energy density of matter fields, motivating the non-negativity assumption.
\end{remark}

General initial data for the Einstein equations is given by a Riemannian manifold $(M,\gamma)$ equipped with a symmetric two-tensor $K$. The time-symmetric condition is the condition that $K$ is identically zero.
\begin{definition}
	We say $(M,\gamma,K)$ is asymptotically flat initial data if $(M,\gamma)$ is an asymptotically flat manifold as given by Definition \ref{def:AF}, and $K_{ij}$ is a symmetric two-tensor satisfying $K=O(|x|^{-\tau-1})$ and $\onabla K=O(|x|^{-\tau-2})$. 
\end{definition}
\begin{definition}
	The ADM mass of an asymptotically flat initial data set $(M,\gamma,K)$ is given by $\m_{ADM}(M,\gamma,K)=E^2-|p|^2$ where $E$ is given exactly by the expression \eqref{eq:admmassdefn} for the ADM mass in the time-symmetric case -- here it is often called the ADM energy -- and $p\in \mathbb R^{n}$ is the ADM (linear) momentum given by
	\begin{equation}
		p_i=\frac{1}{(n-1)\omega_{n-1}} \lim_{R\to\infty}\int_{S_R} (K_{ij}-\tr_\gamma(K)\gamma_{ij})\nu^j\,dV.
	\end{equation}
\end{definition}
\begin{definition}\label{def:constraints}
	The constraint map, $\Phi$ is given by
	\begin{align} \begin{split}
			\Phi_0(\gamma,K)&=R(\gamma)+(\tr_\gamma K)^2-|K|^2\\
			\Phi_i(\gamma,K)&=2( \nabla^jK_{ij}-\nabla_i(\tr_\gamma K)),
		\end{split}
	\end{align}
	and we write $\rho=\Phi_0(\gamma,K)$ and $J_i=\Phi_i(\gamma,K)$. We say initial data $(M,\gamma,K)$ satisfies the \textit{dominant energy condition} (DEC) if $\rho\geq \sqrt{\gamma^{ij}J_i J_j}$. This is equivalent to non-negative scalar curvature in the time-symmetric case.
\end{definition}

\begin{thm}\label{thm:PMTST}
	Let $(M,\gamma,K)$ be complete asymptotically flat initial data (without boundary) satisfying the DEC and with $R(\gamma)\in L^1(M)$, then the ADM mass satisfies $\m_{ADM}(M,\gamma,K)\geq0$ with equality if and only if $(M,\gamma,K)$ corresponds to a hypersurface in $(n+1)$-dimensional Minkowski space.
\end{thm}

\begin{remark}\label{rem:pmt}
	The Positive Mass Theorem was first proven by Schoen and Yau \cite{SY1, SY2}, and independently by Witten \cite{wittenPMT}, both in dimension $3$ and using entirely different techniques. Both proofs have been since generalised to higher dimensions, the former up to $n=7$ \cite{PMTdim8} and the latter for all dimensions provided the manifold is spin \cite{Bartnik-86}. The remaining cases are somewhat more delicate and have since been treated by Lohkamp \cite{lohkamp1, lohkamp2}, and independently by Shoen and Yau \cite{SYPMT2001}. The rigidity statement was also recently addressed for non-spin manifolds by Huang and Lee \cite{huangleepmt1,huangleepmt2}.
\end{remark}

\begin{remark}
	The linearisation of the constraint map plays an analogous role to the linearisation of scalar curvature, in understanding the ADM mass outside of time-symmetry.
\end{remark}

Although the background and definitions in this section are given for arbitrary dimension, we will focus on the special case of $n=3$ for the remainder of this article as most of the work on the topic focuses on this case. Nevertheless, most of the results discussed should also hold in higher dimensions. Some attempt will be made to remark where higher dimensions have been explicitly considered, or where dimension $3$ is used in an essential way.

We take this opportunity to mention that throughout we will predominantly use $\gamma$ to denote a metric on a $3$ manifold $M$, and $g$ to denote a metric on a closed surface $\S$.

\section{Defining the Bartnik mass}\label{sec:defn}

\subsection{A first definition}
We review the various definitions of Bartnik mass and how they relate to each other, as several subtly different variations exist throughout the literature. We begin with the first definition given by Bartnik \cite{Bartnik-89}, and focus on the time-symmetric case. In fact, we will focus predominantly on the time-symmetric case in this review -- not as a matter of preference, but simply because very little is understood outside of time-symmetry. Time-symmetric initial data is simply a Riemannian manifold $(M,\gamma)$, which we usually assume has non-negative scalar curvature.

Given a domain $\Omega$ with connected boundary $\partial \Omega$ in an asymptotically flat manifold $(M,\gamma)$, we define the set $\mathcal A_o(\Omega,\gamma_\Omega)$ of \textit{admissible extensions} to $\Omega$ to be the set of all asymptotically flat manifolds with non-negative scalar curvature in which $(\Omega,\gamma_\Omega)$ isometrically embeds, with no stable closed minimal surfaces outside of the image of $\Omega$. When there is no potential for ambiguity, we will abuse notation slightly to write $\Omega$ to indicate both the original subset of $M$ and its isometric embedding in other manifolds. The first definition of Bartnik mass, $\m_o$, in the time-symmetric case is then
\begin{equation}
\m_1(\Omega,\gamma_\Omega)=\inf\{ \m_{ADM}(M,\gamma):(M,\gamma)\in \mathcal A_1(\Omega,\gamma_\Omega) \},
\end{equation} 
where $\m_{ADM}(M,\gamma)$ is the ADM mass of $(M,\gamma)$. The exclusion of closed minimal surfaces is often referred to as a \textit{non-degeneracy condition}, and while this was the original formulation given by Bartnik, other non-degeneracy conditions are frequently considered. Such a condition is required to exclude extensions wherein $\Omega$ is shielded from infinity by a horizon. Without such a condition, the Bartnik mass of any domain would be zero, as one can construct an extension that contains an arbitrarily small minimal surface enclosing $\Omega$, and then simply glue a small mass Schwarzschild manifold to the minimal surface. We discuss the subtleties of the non-degeneracy conditions later.

\subsection{Definition from Bartnik boundary data}\label{def:boundarydata}
As mentioned in the introduction, we cannot expect that the infimum, if it is realised, be smooth across $\p\Omega$. It is therefore reasonable to expand $\mathcal A_1$ to a larger set of extensions that are only Lipschitz across $\p\Omega$. In this case the scalar curvature must be taken to be non-negative in the distributional sense, which can be understood via the Gauss equation as follows. Consider a small neighbourhood of $\S_0=\p\Omega$ foliated by level sets of the signed distance to $\S_0$, and assume for now that the metric is in fact smooth. By the Gauss equation and the second variation of area formula, the scalar curvature of $(\Omega,\gamma)$ near $\S_0$ is given by
\begin{equation*}
	R(\gamma)=R(\S_t)-|\Pi_t|^2-H_t^2-2\nabla_{\nu_t}(H_t),
\end{equation*}
where $R(\S_t)$, $\Pi_t$ and $H_t$ are the scalar curvature, second fundamental form and mean curvature of $\S_t$ respectively, and $\nu_t$ is the unit normal vector with orientation consistent with $H_t$. From this we can see that if $\gamma$ were only Lipschitz along $\S_0$, then in order to avoid a Dirac delta-type spike in scalar curvature, we must have that the mean curvature on each side of $\S_0$ be equal. Alternatively, we could simply impose that the mean curvature can only decrease across $\S_0$ to ensure the distributional spike contributes positively to the scalar curvature.

That is, we ask that the (outward-pointing) mean curvature of $\partial \Omega$ with respect to the interior and exterior of $\partial \Omega$, which we denote by $H_-$ and $H_+$ respectively, satisfy $H_-\geq H_+$. This allows us to consider extensions of $\Omega$ to be asymptotically flat manifolds $M$ with boundary $\partial M$, such that the induced metric on the boundary agrees with the induced metric $g$ on $\partial \Omega$. One is able to glue such an $M$ to $\Omega$ by identifying $\partial M$ with $\partial \Omega $ to obtain an asymptotically flat Lipschitz manifold $\widehat M$ without boundary (see \cite{miao2002} for details). This leads one to then consider another space of admissible extensions $A_2(\S,g,H)$ to be the set of asymptotically flat manifolds with non-negative scalar curvature, inward-pointing (that is, towards the asymptotic end) boundary mean curvature $H$, and boundary isometric to $(\S,g)$, satisfying a non-degeneracy condition like the one mentioned above. This leads to a definition of (time-symmetric) Bartnik mass in terms of \textit{Bartnik data} $(\S,g,H)$ related to the previous definition by considering the closed manifold $\S$ to be $\partial \Omega$ with $g$ and $H$ the induced metric on $\S=\partial \Omega$ and the outward-pointing mean curvature respectively. Namely, we define the Bartnik mass in terms of Bartnik data by
\begin{equation}
\m_2(\S,g,H)=\inf\{ \m_{ADM}(M,\gamma):(M,\gamma)\in \mathcal A_2(\S,g,H) \}.
\end{equation} 
Note that this is defined for an arbitrary closed manifold $(\S,g)$ equipped with a function $H$, however if this Bartnik data does not come from a domain $\Omega$ in an asymptotically flat manifold then we can not be sure that $\mathcal A_1$ is non-empty, nor can we ensure the Bartnik mass is non-negative. However, this definition has several advantages over the original definition, as it allows us to view the construction of extensions as a kind of boundary value problem. Most notably is that considering extensions in the sense of $\mathcal A_2$ allows one to approach the static metric extension conjecture as a boundary value problem for the static initial data equations. This is discussed in more detail in Section \ref{sec:statstat}.

Note that in the above, the condition that scalar curvature be non-negative distributionally implies that one should only ask that the extensions considered have boundary mean curvature bounded above by the given $H$ from the Bartnik data. That is, we could define $\mathcal A_3(\S,g,H)$ similarly to $\mathcal A_2(\S,g,H)$ except we allow $(M,\gamma)\in \mathcal A_3(\S,g,H)$ to have (inward-pointing) boundary mean curvature $H_+\leq H$. Then $\m_3(\S,g,H)$ is defined in terms of this. This is a less useful definition from the perspective of the boundary value problem, although it is closer in spirit to the original definition. Nevertheless, it is well-understood that positive scalar curvature contributes positively to the ADM mass, so one should expect that $\m_3=\m_2$ generically, and indeed this appears to be the case (see Section \ref{sec:equivs} below).

\subsection{Non-degeneracy conditions}

Although the original definition required that extensions exclude stable closed minimal surfaces, other non-degeneracy conditions often replace this throughout the literature. Generally this is because other conditions can be equally well-motivated, but may exhibit certain useful properties for the problem at hand. For example, one could rule out all minimal surfaces instead of only stable ones, since if an extension contains a minimal surface that is not stable one should be able construct an extension with arbitrarily close ADM mass and without the minimal surface via a perturbation argument. In fact, although the stability property better illustrates the motivation for the non-degeneracy condition, it is often omitted in the literature. We will call the condition that an extension contains no closed minimal surfaces enclosing the Bartnik data (or domain $\Omega$) \textit{non-degeneracy condition (A)}.

Another non-degeneracy condition that often appears in the literature is that $\partial M\cong\partial \Omega$ be \textit{outer-minimising} in the extension. That is, there cannot contain any surface homologous to $\partial M$ with less area. We will refer to this as \textit{non-degeneracy condition (B)}. A key advantage this condition has over condition (A) is that the Hawking mass of a closed area outer-minimising hypersurface bounds the ADM mass by below. This turns out to be very useful in establishing many properties of the Bartnik mass, which we elaborate on in Section \ref{sec:criteria}. For this reason, it appears condition (B) is more common in the modern literature.

For a detailed discussion of all the different non-degeneracy conditions that are used throughout the literature, the reader is referred to Section 5 of \cite{jauregui2019smoothing}. In light of the definitions $\m_1,\m_2,\m_3$ and the different non-degeneracy conditions, it is important to understand when these different definitions agree with each other. This is discussed in more detail in Section \ref{sec:equivs}.

\subsection{The spacetime Bartnik mass}\label{sec:STmass}

As noted above, the Bartnik mass is frequently only considered in the time-symmetric case. In fact, it is has become common in the literature that one simply says ``Bartnik mass'' to refer to the time-symmetric version and explicitly refers to the general case as the \textit{spacetime Bartnik mass}. The spacetime Bartnik mass is again formulated in terms of minimising the ADM mass of admissible extensions, however one wishes now to extend a region $\Omega$ in initial data set $(M_o,\gamma_o,K_o$), rather than just a Riemannian manifold. We define a set of admissible extensions $\widehat {\mathcal A_1}$ now as the set of all asymptotically flat initial data sets $(M,\gamma,K)$ such that $\Omega$ can be isometrically embedded into $(M,\gamma)$ with $K_o$ agreeing with the pullback of $K$ on $\Omega$ under this isometry. We also ask that each $(M,\gamma)\in\widehat {\mathcal A_1}$ satisfies the dominant energy condition and satisfies an appropriate non-degeneracy condition in the sense of the preceding section. It should be remarked however, that the appropriate non-degeneracy condition for this case is far less clear. However, The condition Bartnik suggested in \cite{tsinghua} is that the extensions should contain no apparent horizon, which generalises the notion of minimal surfaces to initial data sets.

The Bartnik mass $\widehat{\m_1}$ of $\Omega$ is then defined as the infimum of the ADM mass taken over $\widehat {\mathcal A_1}$. However, to the best of the author's knowledge there is essentially nothing known about this definition directly. That is, the handful of results that are available focus on the boundary data formulation of the space-time Bartnik mass. That is, analogously to the time-symmetric case, one can consider extensions that fail to be smooth across $\p\Omega$ and then with appropriate boundary conditions consider extensions to be manifolds with boundary instead.

To see what the boundary data should be, we follow a similar reasoning as in the time-symmetric case. Note that this follows closely Bartnik's exposition given in \cite{tsinghua}. In order to do this, we recall the constraint map, $\Phi$, given by Definition \ref{def:constraints} and the dominant energy condition that we would like to impose.

Let $\S_0$ denote a $2$-surface enclosing a domain $\Omega$ in a given initial data set $(\widehat M,\widehat g,\widehat K)$, and we again assume a smooth foliation by level sets of the signed distance to $\S_0$, as above. From the second variation of area and Gauss equation again, we have
\begin{equation}
	16\pi \rho = \Phi_0(g,K) = R(\S_t) - |\Pi_t|^2 - H_t^2 - 2\nabla_ {\nu_t} H_t + (\tr_{\S_t}K)^2 - |K|^2.
\end{equation}
We see that to avoid a distributional spike in $\rho$ we still need only to ensure $H$ matches on each side of $\S_0$. The remaining geometric boundary conditions arise from the momentum constraint. For clarity we drop reference to  $t$ and it should be understood that we continue to work on this foliation. It will be helpful to work with coordinates adapted to the foliation, so let ${\partial_A}$ with $A=1,2$ be a frame on $\S$.

From the momentum constraint, we have
\begin{equation*} 
8\pi J_{\nu}=\nabla^A(K_{A\nu})-\nabla_\nu(\tr_\S K);
\end{equation*} 
a tangential derivative that is bounded, and a normal derivative of $\tr_\S K$. This implies that we must ask that $\tr_\S K$ matches on both sides of $\S$ to avoid a distributional spike in $J_\nu$. The tangential components of momentum constraint give
\begin{equation*} 
8\pi J_A=\nabla_\S^B K_{AB}+K_{\nu B} \Pi^{B}_A+K_{\nu A} H+\nabla_\nu K_{A\nu}-\nabla_A(\tr_\S K)-\nabla^\S_AK_{\nu\nu}.
\end{equation*} 
As above, due to the term $\nabla_\nu(K_{\nu A})$, we must ask that $\omega^\perp_A:=K_{\nu A}$ match on either side of $\S$ to avoid a distributional spike in $J_A$.

This motivates the definition of an admissible extension of $\S$ to be an initial data set $(M,g,K)$ with boundary, such that on $\partial M$ the quantities $(g_{\partial M},H,\omega^\perp_A,\tr_{\partial M} K)$ are prescribed by the corresponding quantities on $\S$ in $(\widehat M,\widehat g,\widehat K)$. In a similar spirit to $\m_3$ above, one could instead impose an inequality on the Bartnik data rather than an equality, however this is generally not considered for the same reasons $\m_3$ is generally not considered.

With that said, the space-time Bartnik mass $\widehat {\m_2}(\S,g,H,\omega^\perp,\tr_\S K)$ is given as the infimum of the ADM mass over the space $\widehat{\mathcal A_2}(\S,g,H,\omega^\perp,\tr_\S K)$ of extensions described above.

\subsection{The equivalence of definitions}\label{sec:equivs}
It is an important question to ask when and if each of the subtly different definitions of Bartnik mass are in fact equivalent to each other. However, even in the time-symmetric situation, it is far from clear that each different definition $\m_1$, $\m_2$ and $\m_3$ with respect to each different non-degeneracy condition are in fact equivalent. In fact, not only is it clear that some Bartnik data $(\S,g,H)$ do not correspond to any $\Omega$ coming from an asymptotically flat manifold with non-negative scalar curvature, we cannot even determine in general if given Bartnik data $(\S,g,H)$ can be realised as the boundary of a compact manifold with non-negative scalar curvature. So to even begin to compare the definitions, we should assume that we have some domain $\Omega$ in an asymptotically flat manifold $(M,\gamma)$ with non-negative scalar curvature inducing some Bartnik data $(\S,g,H)$ on the boundary.

It is clear that, for a given non-degeneracy condition, we have $\m_1(\Omega,\gamma_\Omega)\geq\m_2(\S,g,H)\geq \m_3(\S,g,H)$ as the sets of admissible extensions are nested subsets of each other. However, it is not clear when these inequalities are in fact equalities, or when different non-degeneracy conditions yield the same mass. Equality between $\m_2$ and $\m_3$ seems heuristically most obvious because one should expect that it is always possible to decrease the ADM mass of a manifold with compact boundary by increasing the boundary mean curvature. However even in this case it does not appear to be fully resolved. We do however have the following recent result due to Jauregui \cite{jauregui2019smoothing}.
\begin{thm}[Jauregui]\label{thm:equiv}
	Suppose $\Omega$ is a connected compact subset of an asymptotically flat manifold $(M,\gamma)$ with non-negative scalar curvature, such that the boundary $\partial \Omega$ has positive mean curvature (with respect to the outward-pointing normal). If the non-degeneracy condition is taken to be that $\S \cong \partial \Omega$ is strictly area outer-minimising in the extensions, then all three definitions of Bartnik mass given above yield the same value. That is,
	\begin{equation*}
		\m_1(\Omega,\gamma_\Omega)=\m_2(\S,g,H)=\m_3(\S,g,H),
	\end{equation*}
	where $g$ is the induced metric on $\partial \Omega$.
\end{thm}
\begin{remark}
	Note that the non-degeneracy condition is not exactly condition (B), but rather requires that the outer-minimising condition be \textit{strict}.
	The key property of this choice of non-degeneracy condition is that it is an open condition, in the sense that it is stable to small perturbation. This makes it possible to deform extensions without violating the non-degeneracy condition.
\end{remark}

\begin{remark}
	Independently, equivalence of definitions was shown in \cite{mccormick2020gluing} where the non-degeneracy condition was taken to be the usual condition (B), however additional much stronger convexity conditions were imposed.
\end{remark}

In light of this, in what follows we will regularly use $\m_B$ to denote the Bartnik mass, and only distinguish between the definitions where it is of significance.

It should be remarked though, that outside of time-symmetry, it is still an open question whether the boundary data formulation is in fact equivalent to the definition in terms of isometric embeddings. However, it seems likely to be the case, at least with a reasonable non-degeneracy condition imposed.

\begin{remark}
	Although dimension $3$ is motivated by physics and the focus here, we remark that the Bartnik mass can be easily formulated in higher dimensions, and one expects that most results on the topic can be generalised to higher dimensions. Comments will be made throughout regarding obstructions or complications in the case of higher dimensional analogues.
\end{remark}

\section{Properties of a good quasi-local mass}\label{sec:criteria}

As mentioned above, there is a relationship between the Hawking mass and the ADM mass that one can exploit to establish several properties of the Bartnik mass. Namely, if $(M,\gamma)$ is an asymptotically flat $3$-manifold with non-negative scalar curvature and $\S$ is an outer-minimising surface in $M$, then Huisken and Ilmanen's proof of the Riemannian Penrose inequality \cite{HI} implies $\m_H(\S)\leq \m_{ADM}(M,\gamma)$. Recall the Hawking mass is given by

\begin{equation*}
\m_H(\S)=\left(\frac{|\S|}{16\pi}\right)\left(1-\frac{1}{16\pi}\int_\S H^2\,d\S\right),
\end{equation*} 

where $|\S|$ is the area of $\S$. Assuming non-degeneracy condition (B), this immediately implies that the Bartnik mass is bounded below by the Hawking mass. It is important to note that the proof of this relationship between the Hawking mass and ADM mass (therefore also Bartnik mass) heavily relies on the dimension being $3$, as it employs the Gauss--Bonnet theorem on hypersurfaces. We are now in a position to relate the Bartnik mass to the properties (I) -- (V) given in the introduction.\\

\noindent \textbf{(I) \textit{Positivity}.}\\[-2mm]

Since the Bartnik mass is defined via the ADM mass, which satisfies the positive mass theorem, every admissible extension in the sense discussed above must have non-negative mass if the Bartnik data comes from a domain in an asymptotically flat initial data set satisfying the dominant energy condition. This property is satisfied for the most commonly discussed quasi-local mass definitions with one notable exception. The Hawking mass, which is straightforward to calculate and an incredible tool in geometric analysis, fails this condition spectacularly for surfaces that are far from round spheres.\\

\noindent \textbf{(II) \textit{Compatibility with the ADM mass}.}\\[-2mm]

This was established in the time-symmetric case for the Bartnik mass by Huisken and Ilmanen \cite{HI}, as a consequence of the inequality they obtained between the Hawking mass and the Bartnik mass, mentioned above. The idea is that for any given exhaustion $\{\Omega_i\}$ of a manifold $(M,\gamma)$, for sufficiently large $i$ there is a sequence of spheres $S_i\subset\Omega_i$ with the Hawking mass of $S_i$ converging to the ADM mass. Assuming non-degeneracy condition (B) to ensure that the Bartnik mass is bounded below by the Hawking mass and above by the ADM mass, the conclusion follows. Although the relationship between the Hawking mass and the ADM mass is only known in dimension 3, a lower bound for the ADM mass due to Miao and the author \cite{MM19penrose} (see Theorem \ref{thm:fillins} below) can be used in place of the Hawking mass to establish the large sphere limit in asymptotically Schwarzschild manifolds in dimensions $4\leq n\leq7$, following similar arguments to those used by Jauregui in \cite{jaureguisemicont2}. The question remains open outside of time-symmetry. \\

\noindent \textbf{(III) \textit{Monotonicity}.}\\[-2mm]

This is almost in some sense baked into the definition of the Bartnik mass for free, since if $\Omega_1\subset\Omega_2\subset M$ then every extension of $\Omega_2$ gives an extension of $\Omega_1$. However, there is a subtle problem here coming from the non-degeneracy condition. It can be seen from the definition that in the time-symmetric case, and with non-degeneracy condition (B), then we do have monotonicity under the assumption that $\Omega_1$ is outer-minimising in $\Omega_2$. However, for non-degeneracy condition (A), and outside of time-symmetry, this remains a subtle issue to be considered.\\

\noindent \textbf{(IV) \textit{Small sphere limit}.}\\[-2mm]

There is a good expectation for what quasi-local mass should recover in the limit when calculated on small spheres, such as to highest order recovering the energy density of any matter fields. While there are some subtleties surrounding how to define this small sphere limit, for our purposes we will mean the limit of the quasi-local mass computed on geodesic spheres $S_r$ in an initial data set, of radius $r$ as $r$ tends to zero. We will again only discuss the time-symmetric case, as there are no small sphere results for the Bartnik mass in the more general case. It is quite standard that one expects the small sphere limit of quasi-local mass to give an expression like
\begin{equation*}
	\m_{QL}(S_r)=\frac{4}{3}\pi r^3 T_{00}+O(r^5)=\frac{R}{12}r^3+O(r^5),
\end{equation*}
where $T_{00}=\frac{R}{16\pi}$ is the local energy density of the matter fields and $R$ is the scalar curvature of the metric. The $O(r^5)$ comes from the Bel--Robinson Tensor, which originates from an attempt to construct a local stress-energy density for the gravitational field analogously to the stress-energy tensor of electromagnetism \cite{Bel1959}. Most of the quasi-local mass definitions considered in the literature successfully recover the matter energy density with the $r^3$ term (see, for example, \cite{fanshitam, horowitzschmidt, LRR}), however there is less certainty surrounding the exact form of the $O(r^5)$ term and there appears to be some disagreement in the literature on this. Most definitions do recover something related to the Bel--Robinson tensor, however this term appears to be more sensitive to how the limit is taken.

In the case considered here, the Bartnik mass in time-symmetry, Wiygul was the first to examine this limit \cite{wiygul1}. Using the small sphere limit of the Hawking mass \cite{fanshitam, horowitzschmidt} as a lower bound for the Bartnik mass and by constructing static extensions with controlled mass for an upper bound, he obtains
\begin{equation*}
	\m_B(S_r)=\frac{R}{12}r^3+O(r^4).
\end{equation*}
Note that Wiygul considers non-degeneracy condition (B) and takes the mass to be $\m_2$, however it appears that the non-degeneracy condition could easily be replaced by a strictly outer-minimising condition so that Jauregui's work on the equivalence of definitions would then apply here (Theorem \ref{thm:equiv}).

In a subsequent paper Wiygul \cite{wiygul2} improved the upper bound to better control the higher order terms, which combined with a lower bound from the known small sphere limit on the Hawking mass gives
\begin{align}\label{eq:smallsphere}\begin{split}
	\left(\frac1{120}\Delta R-\frac{R^2}{144}  \right)r^5+O&(r^6)\leq \m_B(S_r)-\frac{R}{12}r^3\\
	&\leq
	\left(\frac1{120}\Delta R-\frac{5R^2}{432}+\frac1{72}|\ric|^2  \right)r^5O(r^6),\end{split}
\end{align}
where $\ric$ is the Ricci curvature. Harvie and Wang \cite{harviewangspheres} have since confirmed this upper bound using different techniques. Since the upper bound comes from a static extension, which one expects to realise the Bartnik mass in some generality (see Section \ref{sec:statstat}), is seems likely that the upper bound from \eqref{eq:smallsphere} is in fact an equality. In fact, it also seems reasonable that the small sphere limit for the Hawking mass could in fact agree with this upper bound when a different limiting process is used. 

The Hawking mass lower bound again presents an issue for generalising this to higher dimensions, however it is likely that the analogous inequality obtained in \cite{MM19penrose} can again be used in place of the Hawking mass lower bound, under certain hypotheses on the point around which the spheres are centred.\\

\noindent \textbf{(V) \textit{Rigidity}.}\\[-2mm]

In the time-symmetric case, this property also follows from the work of Huisken and Ilmanen \cite{HI}, both relating the Hawking mass to the ADM mass and directly making use of inverse mean curvature flow. As above, this requires that one uses non-degeneracy condition (B) in the definition. In particular, Huisken and Ilmanen show that if some domain $\Omega$ contains any point $x$ where the manifold is not flat, there exists an inverse mean curvature flow from $x$ out to infinity such that the Hawking mass is strictly positive along the flow. Then by showing there exists a leaf of the flow close to $x$ that is outer-minimising, positivity of the mass in any extension follows.

Interestingly, the result is not that the domain $\Omega$ must be a domain in Euclidean space, only that it is \textit{locally} isometric to it. It turns out, by the work of Anderson and Jauregui \cite{A-J}, this can be shown to be the best one can hope for. Specifically, they show the existence of a domain with zero Bartnik mass such that it cannot be isometrically embedded in $\mathbb R^3$. Note that by the rigidity of the positive mass theorem (with corners \cite{miao2002}), this implies that the infimum is not realised by any admissible extension. We will revisit this point below.

\section{Stationary and static extensions}\label{sec:statstat}
Although the Bartnik mass seems impossible to compute in general, there is some hope that this will not always be the case. In particular, if we knew the mass was realised by a particular extension that can be uniquely determined from the boundary data, and it's mass computed from there, we would be able to compute the Bartnik mass directly.

In this direction, Bartnik conjectured \cite{tsinghua} that, at least under reasonable hypotheses, that the Bartnik mass is realised as the ADM mass of a stationary extension to the domain/Bartnik data. That is, initial data that generates a spacetime with a Killing field that is timelike at infinity. Again we begin with the time-symmetric case, which is more well-understood. There the conjecture is that the metric extension be static rather than stationary, which one can formulate at the initial data level as follows. An asymptotically flat Riemannian manifold $(M,\gamma)$ is said to be \textit{static} if there exists a positive solution $N$ to the system
\begin{equation}\label{eq:static}
	\begin{split}
		\Delta N&=0\\
		\nabla_i\nabla_j N&=N\ric(\gamma)_{ij}
	\end{split}
\end{equation}
on $(M,\gamma)$. In this case, we say $N$ is a \textit{static potential}. Note that if a solution $N$ does vanish somewhere, then the zero set is a totally geodesic hypersurface (\cite[Prop. 2.6]{corvino2000scalar}). Additionally, it is often required that $N$ be asymptotically constant and scaled to equal $1$ at infinity for physical reasons. Note that by work of Beig and Chru\'sciel \cite{beigchrusciel} (see also, Miao and Tam \cite[Prop. 3.1]{miaotam}) a bounded static potential on an asymptotically flat manifold is always asymptotically constant.

Partially due to the subtleties of different mass definitions, and partially because it is unclear what conditions the conjecture should require, we avoid stating the conjecture here in a more precise form and simply summarise some partial results that are known. In fact, The most optimistic version of this conjecture is already known to be false by the work of Anderson and Jauregui, mentioned in the preceding section. There, they demonstrate the existence of domains with zero Bartnik mass where the infimum is not realised \cite{A-J} (see also the very recent work on this by Anderson \cite{Anderson24}). In light of this example it seems reasonable to believe one could construct similar examples of domains with non-zero Bartnik mass coming from non-embedded hypersurfaces in static manifolds, however to the best of the author's knowledge this remains a problem to be checked.

\subsection{Static metrics as mass minimisers}
In the time-symmetric case, where the conjecture is that the mass be minimised by a static metric, the problem is quite well understood and closely related to scalar curvature deformation. To see this, recall the linearised scalar curvature given by \eqref{eq:DR}. One can check that, on an asymptotically flat manifold, \eqref{eq:static} is equivalent to the equation $D_\gamma R^*(N)=0$, where $D_\gamma R^*$ is the formal adjoint of $D_\gamma R$. The implication \eqref{eq:static} $\implies D_\gamma R^*(N)=0$ is straightforward to check, however the other implication first requires showing that $D_\gamma R^*(N)=0$ implies that $R$ must be constant \cite[Theorem 1]{FM1975deformations} (see also \cite[Proposition 2.3]{corvino2000scalar}), then asymptotic flatness implies that constant must be zero. It was shown by Corvino \cite{corvino2000scalar} that an asymptotically flat manifold that is not static can be deformed locally increasing scalar curvature which after a conformal change to bring it back down leads to a decrease in the ADM mass. Later he showed that this can be done on a manifold with boundary while preserving a neighbourhood of the boundary \cite{corvino2017note}, which therefore proves that an admissible asymptotically flat extension that is not static can always have its mass reduced by a small amount. This was also established independently and by different methods by Anderson and Jauregui \cite{A-J}. The precise statement given by Corvino is the following.
\begin{thm}
	Suppose $(M,\gamma)$ is a smooth connected Riemannian $n$-manifold with compact boundary $\S = \partial M$, and the scalar curvature $R(\gamma) \in L^1(M)$ is non-negative. Further suppose $D_\gamma R^*$ has trivial kernel on $M$.
	
	Then there is an $\varepsilon_0 > 0$, an open set $\Omega$ containing $\S$, and a family ${\gamma_\varepsilon}$, with $\varepsilon\in(0,\varepsilon_0)$, of asymptotically flat metrics on $M$ satisfying:
	\begin{itemize}
		\item $\gamma_\varepsilon = \gamma$ on $\Omega$,
		\item $R(\gamma_\varepsilon) = R(\gamma)$,
		\item as $\varepsilon \to 0$, $\gamma_\varepsilon$ converges $\gamma$ in $C^k$ on compact subsets (for any positive integer $k$), and
		\item  the ADM mass of $\gamma_\varepsilon$ is strictly less than the ADM mass of $\gamma$.
	\end{itemize}
\end{thm}
Note that this result holds in all dimensions $n\geq3$. As noted by Anderson and Jauregui \cite{A-J}, this in itself is not quite sufficient to prove that a mass minimising extension must be static, if it exists. This is because one must take care to ensure the non-degeneracy condition is preserved. However, if one assumes a non-degeneracy condition that is stable to perturbations -- for example, the strict version of condition (B) discussed above -- then this suffices to prove that mass-minimising extensions must be static.

As mentioned above, independently of Corvino's work, Anderson and Jauregui \cite{A-J} established this by other methods. Essentially they completed the approach put forward initially by Bartnik \cite{bartnik2005}, which equates critical points of the ADM mass along a manifold of scalar flat metrics with static solutions (see also \cite{mccormick2015note,mccormick2021hilbert}). Morally the argument is a Lagrange multipliers argument based on the functional $16\pi\m_{ADM}-\int_M NR(\gamma)\,dV_\gamma$, where $N$ plays the role of a Lagrange multiplier yielding a static potential. They prove that an extension of given Bartnik data $(\S,g,H)$ that is a critical point of the ADM mass, over a space of admissible extensions modelled on weighted H\"older spaces, must be static. Furthermore, that the converse is true. In higher dimensions, some work in this direction was recently carried out by Delay \cite{delay2020bartnik}.

\subsection{Stationary metrics as mass minimisers}
More recently it has been shown that both of the approaches mentioned above can be adapted to prove that outside of time-symmetry, a mass-minimising vacuum extension must be stationary. To better understand this, we must reformulate the constraint map in terms of the momentum $\pi^{ij}=K^{ij}-\tr_\gamma(K)\gamma^{ij}$. That is, we consider the map
\begin{equation*}
	\widehat\Phi(\gamma,\pi)=\begin{bmatrix}
		\Phi_0(\gamma,K)\\ \Phi_i(\gamma,K)
	\end{bmatrix}.
\end{equation*}
It is now well-known that vacuum initial data $(\gamma,\pi)$ for which there exists a non-trivial solution $(N,X^i)$ to
\begin{equation}\label{eq:stationary}
	D_{(\gamma,\pi)}\widehat\Phi^*(N,X^i)=0
\end{equation}
generates a stationary solution to the Einstein equations \cite{moncrief1975}. Corvino showed \cite{corvino2019short}, analogous to the time-symmetric case, that if \eqref{eq:stationary} has no non-trivial solutions then the mass of an extension can be perturbed smaller. That is, again assuming the non-degeneracy condition is stable to perturbations, a mass-minimising extension must admit a non-trivial solution to \eqref{eq:stationary}. This follows from his joint work with Huang \cite{corvino2020localized} where they show that the non-existence of solutions to \eqref{eq:stationary} allows for local deformations of the initial data so that the dominant energy condition is strictly satisfied in a region analogous to Corvino's earlier work on scalar curvature deformation \cite{corvino2000scalar}. Once strictness in the dominant energy condition is achieved, this is used to push the mass down by a small amount.

Independently, An \cite{an2020} established the same result using a similar strategy to Anderson and Jauregui in the time-symmetric case. Namely that if the Bartnik mass is realised by an admissible vacuum extension then that extension must be stationary.

\begin{remark}
	In this case, there is a subtlety that is easy to overlook. Namely the appearance of the word ``vacuum'' in the above discussion. In the time-symmetric case, the existence of a non-trivial solution to \eqref{eq:static} implies vanishing scalar curvature, that is, vacuum. For non-trivial solutions to \eqref{eq:stationary}, no such result is available. However, recently Huang and Lee \cite{huanglee} proved that an extension realising the Bartnik mass must satisfy several related properties. In particular, it corresponds to an initial data slice of a null dust spacetime that admits a global Killing field. 
\end{remark}

\subsection{Existence and uniqueness of static and stationary extensions}\label{sec:existunique}
Given that the Bartnik mass -- if realised by an admissible extension -- must be static (resp. stationary), it is important to understand when such static (resp. stationary) extensions exist, and if they are unique. An important feature of the boundary value definition of Bartnik mass is that is that the problem of existence and uniqueness of static (resp. stationary) extensions becomes a boundary value problem. That is, in the time-symmetric case, we can view \eqref{eq:static} as a boundary value problem for $(\gamma,N)$ on some manifold $M$, where the Bartnik data $(\S,g,H)$ provides the boundary conditions, and $g$ is an asymptotically flat metric inducing $g$ and $H$ on $\partial M= \S$. Outside of time symmetry, the problem would then be viewing \eqref{eq:stationary} as a boundary value problem for $(\gamma,K,N,X)$ where the spacetime Bartnik data $(\S,g,H,\omega^\perp,\tr_\S K)$ is prescribed. One then would like to know, when are these boundary value problems actually well-posed? We first discuss the time-symmetric case. Work of Miao \cite{miao2003existence} shows the existence of static extensions for Bartnik data sufficiently close to that of a sphere in $\mathbb R^3$ under an additional symmetry assumption. Anderson and Khuri \cite{andersonkhuri} later proved that the boundary value problem is elliptic and of Fredholm index zero. Furthermore, Anderson \cite{anderson2015} proved existence and uniqueness of static extensions for Bartnik data close to that of the round sphere in $\R^3$ without the symmetry assumption of Miao. More recently, An and Huang \cite{anhuang} proved existence and local uniqueness for the Bartnik data close to that of more general data coming from $\R^3$. The precise statement of the result requires some additional terminology so the reader is directed to \cite{anhuang} for more details. Uniqueness of static extensions to round CMC data under an additional strict stability-type condition was very recently addressed by Harvie and Wang \cite{HW}. Outside of time-symmetry, as usual quite little is known. However, recent work of An \cite{an2020ellipticity} has established ellipticity of the boundary value problem.

The reader is directed to the survey of recent progress on the Bartnik mass given by Anderson \cite{anderson2020recent}, wherein this approach is discussed in more detail.

\section{Estimates and constructions of extensions}\label{sec:estimates}

Under non-degeneracy condition (B), the inequality between the Hawking mass and the Bartnik mass is the main lower bound available. Furthermore, for round spheres in $\R^3$ or Schwarzschild manifolds, the Hawking mass is equal to the ADM mass, so it therefore is also an upper bound for the Bartnik mass (by virtue of the fact that these provide admissible extensions).

On the other hand, as the Hawking mass is negative for surfaces in $\R^3$ that are not round -- and more negative the further from round the surface is -- one expects this inequality to be strict in general. Once again in higher dimensions where the relationship between the Hawking mass and Bartnik mass is not available, the inequality established in \cite{MM19penrose} can be used instead to conclude equality of the Hawking and Bartnik mass for spheres in Schwarzschild manifolds.

Since the Bartnik mass is taken to be the infimum of the ADM mass over a space of asymptotically flat extensions, one can estimate the Bartnik mass from above by simply constructing one such extension with known ADM mass. To this end, Mantoulidis and Schoen \cite{M-S} constructed extensions of Bartnik data corresponding to stable minimal surfaces, whose mass is arbitrarily close to the lower bound given by the Riemannian Penrose inequality. Their result is stated below including a slight generalisation due to Chau and Martens \cite{chaumartens}, who relax the hypothesis on the eigenvalue from strict positivity to non-negativity.
\begin{thm}[Mantoulidis--Schoen, Chau--Martens]
	Let $g$ be a metric on the $2$-sphere $\mathbb S^2$ such that the first eigenvalue of $-\Delta_g+K(g)$ is non-negative, where $K(g)$ is the Gauss curvature of $g$. Then for any $\varepsilon>0$ there exists an asymptotically flat manifold $(M,\gamma)$ with minimal surface boundary isometric to $(\mathbb S^2,g)$ and ADM mass $m=\m_{ADM}(g)<\left(\frac{A}{16\pi}\right)+\varepsilon$, where $A$ is the area of $(\mathbb S^2,g)$. Furthermore, $(M,\gamma)$ is foliated by mean convex spheres and isometric to a mass $m$ Schwarzschild manifold outside a compact set.
\end{thm}
Dimension plays a role here in two ways. First, the lower bound from the Riemannian Penrose inequality is only established up to dimension 7. More seriously though, the construction requires a path of metrics on a sphere connecting the given metric to a round metric, satisfying certain properties. And while such a path is always possible to construct on $2$ dimensional spheres (for use in constructing $3$-dimensional extensions), in higher dimensions this does not appear to be the case. This higher dimensional version was considered explicitly by Cabrera Pacheco and Miao \cite{CM}, who give details on when this construction can be adapted to higher dimensions.

For Bartnik data with non-zero mean curvature, extensions have also been constructed to provide upper bounds on the Bartnik mass. There are essentially two different approaches that have been used recently to construct such extensions. One is an adaptation of the Mantoulidis and Schoen work mentioned above, while the other employs Ricci flow to construct the extensions. As all of the estimates are similar in qualitative properties, we avoid stating each of them precisely and instead refer the reader to relevant references. The Ricci flow approach was used by Lin and Sormani \cite{L-S}, while the approach following \cite{M-S} was first used by Cabrera Pacheco, Cederbaum, Miao and the author \cite{CCMM}. The latter was also refined by several authors \cite{chaumartens2,miaopiubello,miaowangxie} (see \cite{CCsurvey} for an overview and survey of the different applications of this technique). In both approaches, the Bartnik data was assumed to have constant mean curvature (CMC) and the upper bound approaches the Hawking mass as the Bartnik data becomes closer to round in some sense. Both approaches also share the unfortunate feature that the upper bounds are not explicitly computable from the Bartnik data itself -- it remains an interesting problem to obtain more explicit upper bounds on the Bartnik mass. Note that when the definitions $\m_2$ and $\m_3$ of Bartnik mass are indeed equivalent (see Section \ref{sec:equivs}), estimates for CMC Bartnik data provide estimates for non-CMC Bartnik data simply by constructing CMC extensions of $(\S,H_{CMC})$ with $H_{CMC}=\min_\S(H)$ for non-constant $H$.

It should be remarked that an early approach to constructing admissible extensions for the Bartnik mass is the quasi-spherical construction due to Bartnik himself \cite{bartnikquasispherical}. This construction is the basis for the later work of Shi and Tam \cite{shitam} establishing the positivity of the Brown--York mass. In fact, Shi and Tam construct asymptotically flat extensions with prescribed Bartnik data on the boundary and prove that the Brown--York quasi-local mass is bounded below by the ADM mass of the extension. From this, it follows that the Brown--York mass provides an upper bound for the Bartnik mass.

\section{Related quasi-local mass quantities}

\subsection{Bray's inner mass}
In \cite{bray2001} (see also \cite{braychrusciel}), Bray remarks that the Riemannian Penrose inequality suggests a new definition of quasi-local mass in a very similar spirit to the Bartnik mass. In the time-symmetric picture, the Bartnik--Bray inner mass of given Bartnik data $(\S,g,H)$ is defined as follows. Consider a set of admissible fill-ins, which we take to be asymptotically flat manifolds with non-negative scalar curvature and boundary isometric to $(\S,g)$ with \textit{outward} mean curvature $H$. Note, that $H$ taken to be the outward mean curvature -- instead of the mean curvature in the direction of the asymptotic end -- is the reason this is a fill-in, rather than an extension. The inner mass $\m_{inner}$ of given Bartnik data is then taken to be the supremum of the quantity $\frac12(\frac{A}{4\pi})^{1/2}$ over this space of admissible fill-ins, where $A$ is area of the surface of least area in the fill-in that is homologous to $\S$. Outside of time-symmetry, the definition is similar, however we omit a precise definition here, as it does not appear to have been examined in the literature in any detail.

When the standard time-symmetric Bartnik mass is defined with non-degeneracy condition (B), then we have from the Riemannian Penrose inequality that
\begin{equation}\label{eq:innerouter}
\m_{inner}(\S)\leq \m_B(\S).
\end{equation} 
For round spheres in a Schwarzschild manifold, equality in \eqref{eq:innerouter} holds where both quantities are equal to the ADM mass of the Schwarzschild manifold. One can also obtain a strict inequality in \eqref{eq:innerouter} by taking Bartnik data corresponding to a closed surface $\S$ in a Schwarzschild manifold that does not enclose the horizon. In this case, there cannot exist a fill-in with another asymptotically flat end as this would violate the Riemannian Penrose inequality, so $\m_{inner}=0$. On the other hand, the $\m_B>0$ because by \cite{HI} if it were zero the interior must be locally flat, and of course an $m\neq0$ Schwarzschild manifold is not. The reader is directed to Jauregui's work on fill-ins \cite{jauregui2013fill} for some discussion of the inner mass in relation to the construction of fill-ins.

For essentially the same reason that an extension realising the Bartnik mass must be stationary (or static in the time-symmetric case), a fill-in realising the inner mass also must be stationary (or static, respectively). This follows from the work of Chru\'sciel, Isenberg, and Pollack \cite[Theorem 5.1]{CID05}, who therein attributed the observation to Bartnik.

In the same way that the Bartnik mass can be estimated from above by constructing examples of asymptotically flat extensions with controlled mass, the inner mass can be estimated by below by constructing appropriate fill-ins. Motivated by this, Miao and the author proved a ``Penrose-like'' like equality by way of first constructing such a fill-in \cite{MM19penrose}. Although the fill-ins constructed there contain a minimal surface boundary rather than another asymptotically flat end, such an end can be easily constructed by a reflection, giving the following lower bound on the inner mass.
\begin{thm}\label{thm:fillins}
	Let $(\S,g,H)$ be given Bartnik data satisfying $4\min_\Sigma K_g   >  \max_\Sigma H^2$,
	where $ K_g$ is  the Gauss curvature of $(\Sigma, g)$. 
	Then
	\begin{equation}\label{eq-Penrose-like-1}
		\mathfrak{m}_{inner} \geq \frac12 \left(\frac{| \Sigma |}{4\pi}\right)^{\frac{1}{2}} \left( 1 - 
		\frac{  \max_\Sigma H^2 }{ 4\min_\Sigma K_g } \right) , 
	\end{equation}
	where $ | \Sigma | $ is the area of $ \Sigma $.
\end{thm}
Although we state the result in dimension $3$ here, this and another closely related estimate was established for dimensions up to and included $7$.

\begin{remark}
	Sometimes in the literature the definition of the inner mass the definition simply asks that the fill-ins be manifolds with a second boundary component corresponding to a minimal surface, and ask that this be the only closed minimal surface in the fill-in. Then take $A$ to be the area of this minimal surface boundary component, with the same reflection as above in mind.
	
	It was noted by Wang \cite{wang20outerinner} that this definition is essentially the same as the definition of ``outer entropy'' that appears in the AdS/CFT literature. Although in this case, one would consider an asymptotically hyperbolic analogue of the quantity (an analogue of Theorem \ref{thm:fillins} for this case is given in \cite{mccormick2023fillin}). 
\end{remark}

\subsection{Asymptotically hyperbolic analogue}
When the Einstein equations include a cosmological constant $\Lambda$, the Hamiltonian constraint is simply augmented to
\begin{equation*}
	\Phi_o(\gamma,K)=16\pi\rho + 2\Lambda,
\end{equation*}
and the model initial data for isolated systems is hyperbolic space. In this setting, one considers asymptotically hyperbolic manifolds instead of asymptotically flat manifolds, where there exists another appropriate notion of total mass in place of the ADM mass \cite{chruscielherzlich, wang2001mass}. It should be remarked that the mass of such manifolds is somewhat more subtle than the ADM mass, so for the save of exposition we avoid giving a precise definition here. Nevertheless, one can define a notion of mass and then can define a Bartnik mass analogously as the usual definition. That is, the Bartnik mass of a given domain or Bartnik data is taken as the infimum of a hyperbolic total mass over a set of admissible asymptotically hyperbolic extensions. Although the definition is a natural one, it has only recently been considered in the literature.

Some of the known results about the usual Bartnik mass have recently been established for this asymptotically hyperbolic analogue, and it is likely that many other results also have natural analogues to be explored. The reader is directed to \cite{CCM,huangjang,huangjangmartin,martin2023,miaowangxie} for some results on this quantity.

\subsection{A mass for the Einstein--Maxwell system}
Many problems pertaining to mass in general relativity have also been studied in the coupled Einstein--Maxwell system, which includes electromagnetic fields additionally. In this case initial data is supplemented with vector fields $E$ and $B$, representing the electric and magnetic fields, and the constraint equations are supplemented accordingly, including the addition of the Gauss constraint. We avoid too much of a digression into the Einstein--Maxwell system, however make a few remarks in the case of the Bartnik mass in this case.

Given a domain $\Omega$ (or appropriate Bartnik data) in an initial data set for the Einstein--Maxwell equations, one can approach the Bartnik mass in two ways. On one hand, the electromagnetic fields should play no role in the construction, so one can still define the Bartnik mass in the usual way. On the other hand, doing so would mean introducing discontinuities in the electric and magnetic fields that correspond to initial data containing a thin massless shell of electric (or magnetic) charge, which cancel out any electromagnetic fields (and therefore electromagnetic field energy) in the extension. With an alternate ``charged Bartnik mass'' in mind, Alaee, Cabrera Pacheco and Cederbaum \cite{ACC} study the problem of constructing extensions with controlled mass similar to \cite{M-S,CCMM}. A precise formulation of such a charged Bartnik mass was not required there, however a straightforward definition would follow along the lines of Section \ref{def:boundarydata} smoothing out the electromagnetic contributions in such a way to avoid these charged shells -- see, for example, \cite{chenmccormick}, where another ``charged'' quasi-local mass was considered. Although it is not clear that such charged quasi-local mass quantities are of physical significance, it was shown in \cite{chenmccormick} that such a charged quasi-local mass is the correct quantity for a quasi-local Penrose inequality and therefore at least these quantities are interesting from a geometric perspective.

\section{Open problems}
Throughout this article, several open questions have been mentioned. We conclude by presenting a selection of them here.

\subsection{Equivalence of definitions}
As mentioned above, it has been established in the time-symmetric case that under a suitable non-degeneracy condition, the boundary value formulation of the Bartnik mass is equivalent to the definition involving smooth extensions. However, it remains to be understood if (or how) the different non-degeneracy conditions alter the definition. In particular, is the outer-minimising condition equivalent to the \textit{strictly} outer-minimising condition? And are these equivalent to the condition that the extensions contain no closed minimal surfaces enclosing the Bartnik data?

In addition to this, it remains entirely open to understand when and if any of the subtly different definitions agree with each other outside of time-symmetry. It would be of particular value to determine that -- at least under some good choice of non-degeneracy condition -- the boundary data formulation of Bartnik mass agrees with the smooth extension formulation outside of time-symmetry.

\subsection{Estimates outside of time-symmetry}
Following on from above, there is still a considerable amount to understand about the Bartnik mass outside of time-symmetry. No attempt is made to list all of open problems outside of time-symmetry, however one direction that seem interesting is to establish some estimates for the Bartnik mass in this case. For example:
\begin{itemize}
	\item Is there an analogue of the lower bound given by the Hawking mass in time-symmetry? It seems natural to expect that the spacetime Hawking mass should give a similar lower bound, however in light of the fact that Penrose inequality outside of time-symmetry remains an open problem, such a lower bound on the Bartnik mass feels currently far out of reach. Nevertheless, even a crude lower bound would be valuable.
	\item Related to the above, perhaps more optimistically one could hope to find such a lower bound in the case of Bartnik data coming from surfaces far out in the asymptotic end of an asymptotically flat initial data set. Such an estimate would be valuable in proving the large sphere limit of the Bartnik mass outside of time-symmetry -- that is, compatibility with the ADM mass -- which surprisingly remains open.
	\item Can the construction of Mantoulidis and Schoen \cite{M-S} giving an upper bound for the Bartnik mass of a horizon be extended outside of time-symmetry? That is, given Bartnik data corresponding to a general apparent horizon, can an extension be constructed with ADM mass arbitrarily close to the conjectured lower bound coming from the Riemannian Penrose inequality?
	\item Following \cite{CCMM}, it should be possible to construct some crude estimates for the Bartnik mass outside of time-symmetry. However, at the time of writing this no such estimates have been obtained. It would be interesting to find extensions that give estimates for the Bartnik mass tending to the expected value in limiting cases.
\end{itemize}

\subsection{Counterexamples to realising the infimum}
In \cite{A-J}, an example of a domain whose Bartnik mass is not realised by any admissible extension was given. This example has zero Bartnik mass, so one may naturally wonder if this is a unique property of the zero mass case -- that is, an exceptional case -- or if similar counterexamples exist when the mass is non-zero. It does not seem unreasonable that a similar construction could be performed for deformations of domains in a Schwarzschild manifold to construct domains with Bartnik mass of any positive $m$.

\subsection{Explicit upper bounds}
As discussed in Section \ref{sec:estimates} the known upper bounds for the Bartnik mass cannot currently be explicitly computed in terms of the Bartnik data. It would be useful to obtain upper bounds for the Bartnik mass that are explicitly given in terms of the Bartnik data and reduce to the known values for round spheres and minimal surfaces. The estimates outlined in Section \ref{sec:estimates} all involve quantities that are defined abstractly, so one method to obtain explicit estimates would be to estimate these abstract quantities concretely.

\end{document}